\documentclass[11pt]{article}
\usepackage{amsmath,amssymb,amsthm,amsfonts, graphics,color}
\usepackage{psfig}
\usepackage{graphicx,psfrag}
\newcommand{\R}{\mathbb R}
\newcommand{\C}{\mathbb C}

\newcommand\Lh{{\cal L}(H)}
\newcommand\ds{\displaystyle}

\renewcommand{\Re}{\mathop {\rm Re}\nolimits}
\renewcommand{\Im}{\mathop {\rm Im}\nolimits}

\numberwithin{equation}{section}
\newtheorem{theorem}{Theorem}[section]
\newtheorem{ccorollary}[theorem]{Corollary}
\newtheorem{llemma}[theorem]{Lemma}
\begin{document}

\title{Convex domains and $K$-spectral sets}

\author{Catalin Badea, Michel Crouzeix, Bernard
  Delyon}
\date{}
%
\maketitle
%

\begin{abstract}
Let $ \Omega$ be an open convex domain of $\C$. We study
constants $K$ such that $ \Omega$ is $K$-spectral or
complete $K$-spectral for each continuous linear Hilbert space
operator with numerical range included in
$\Omega$. Several approaches are discussed.
\end{abstract}
\section{Introduction}\label{intro}
\subsection{Prologue}
Let $H$ be a complex Hilbert space and let $\Lh$ denote the C$^*$-algebra
of all continuous linear operators on $H$.
For $A\in\Lh$ its numerical range
$W(A)$ is defined by
\begin{eqnarray*}
W(A)=\{\langle Ax,x\rangle \,; x\in H,\,\Vert x\Vert = 1\}.
\end{eqnarray*}
Here $\langle x,y\rangle$ is the inner product in $H$ and $\Vert x\Vert  = \langle x,x\rangle^{1/2}$
the corresponding norm.
Recall that $W(A)$ is a convex subset of $\C$ (the Toeplitz-Hausdorff theorem) and that the spectrum of $A$ is contained in $\overline{W(A)}$. If a rational function $r$ is bounded on the numerical range, then it has no pole in $\overline {W(A)}$; consequently $r(A)$ is
well-defined and belongs to $\Lh$.

\subsection{The constants $C(\Omega)$ and $C_{cb}(\Omega)$}
The aim of this paper is to study the following constant
\begin{eqnarray}\label{co}
C(\Omega):=\sup\{\Vert r(A)\Vert \,; \overline{W(A)} \subset\Omega ,|r(z)|\leq 1,\forall z\in \Omega \},
\end{eqnarray}
where $\Omega$ is an open convex  subset of the complex plane ($\emptyset\neq\Omega\neq\C$).
In this definition the supremum is taken over all complex Hilbert spaces $H$, all continuous linear operators $A$ on $H$ and all rational functions $r\,:\,\C\to\C$ satisfying the prescribed constraints.
Recall also that the operator norm $\|r(A)\|$ is defined by
\begin{eqnarray*}
\|r(A)\|=\sup\{\|r(A)x\|\,; x\in H,\|x\|\leq1\}.
\end{eqnarray*}
In other words $C(\Omega)$ is the best contant $K$ such that the inequality
\begin{eqnarray}\label{spec}
\Vert r(A)\Vert\leq K\, \sup_{z\in \Omega }|r(z)|
\end{eqnarray}
holds for all complex Hilbert spaces $H$, all continuous linear operators $A$
on $H$ with $\overline{W(A)}\subset\Omega$
and all rational functions $r$ bounded in $\Omega$.

It is also interesting to consider the completely bounded analogue of $C(\Omega)$ defined by
\begin{eqnarray}\label{cb}
C_{cb}(\Omega):=\sup\{\Vert R(A)\Vert\, ; \overline{W(A)} \subset \Omega ,\Vert R(z)\Vert \leq 1,\forall z\in\Omega \}.
\end{eqnarray}
As in (\ref{co}), the supremum is taken over all complex Hilbert
spaces $H$, all continuous linear operators $A$ on $H$ (satisfying the prescribed constraints), but now
$R$ runs among the rational functions with matrix values of size $m\times n$ ($R\,:\,\C\to\C^{m,n}$)
and the supremum is taken also over all $m$ and $n$. We refer to the book
\cite{book} for more explanations about completely bounded maps.

We clearly have $1\leq C(\Omega)\leq C_{cb}(\Omega)$ and these constants only
depend on the shape of $\Omega$\,; indeed, they are
invariant under similarities and symmetries. It is also easily seen that the condition
$\overline{W(A)} \subset\Omega$ can be relaxed and that we have
\begin{eqnarray*}
C(\Omega)=\sup\{\Vert r(A)\Vert \,; W(A) \subset\overline\Omega ,|r(z)|\leq 1,\forall z\in \Omega \}\,;
\end{eqnarray*}
a similar formula holds for $C_{cb}(\Omega)$.

Recall that $\Omega$ is said to be a $K$-\emph{spectral set} \cite{book}
(a \emph{spectral set} if $K=1$) for an operator $A$ if the inequality (\ref{spec})
holds for all rational functions without pole in the spectrum of $A$. Therefore our definition (\ref{co}) means
 that $\Omega $\emph{ is a} $C(\Omega)-$\emph{spectral set}
for all operators $A$ such that $\overline{W(A)}\subset \Omega $. Similarly
(\ref{cb}) means that $\Omega $\emph{ is a complete} $C_{cb}(\Omega)-$\emph{spectral set}
for these operators.

A famous result, due to J. von Neumann \cite{vn}, states the following~:

{\it if $\Pi$ is a half-plane, then $\Pi$ is a spectral set for $A$, for any operator $A$ with $W(A)\subset \Pi$;}

\noindent in other words $C(\Pi)=1$. We also have $C_{cb}(\Pi)=1$. These are the first known estimates for the constants $C(\Omega)$ and $C_{cb}(\Omega)$. Since then it has been shown that
\begin{itemize}
\item $C(\Omega)<+\infty$, \qquad for all bounded open convex subsets $\Omega$, \cite{dede},
\item $C_{cb}(S)<2+2/\sqrt3$, \quad for all convex sectors or all strips $S$, \cite{crde},
\item $C(D)=2$, \qquad for all disks $D$, \cite{crzx,okan},
\item $C_{cb}(P)<4.75$, \qquad for all parabolic domains $P$, \cite{cr},
\item $C_{cb}(\Omega)<57$, \qquad for all open convex subsets $\Omega$, \cite{crac}.
\end{itemize}

\noindent Excepting $C_{cb}(\Pi)=1$ and $C(D)=2$, these estimates are far
from being optimal\,; the conjecture $\sup_{\Omega}C(\Omega)=2$
has been proposed in \cite{crzx}, and then the conjecture $\sup_{\Omega}C_{cb}(\Omega)=2$ in \cite{crac}.\smallskip

\subsection{Motivation}
The boundedness of the constants $C(\Omega)$ and $C_{cb}(\Omega)$
allows to extend the rational functional calculus for operators
$A$ satisfying $W(A)\subset \overline\Omega $ to more general
(holomorphic) functions. Furthermore, if $\Omega$ is unbounded,
and after adding a technical resolvent condition, a suitable
functional calculus can be constructed for unbounded operators\,;
we refer for that to \cite{haase}. Up to now, the boundedness of
our constants has allowed to obtain some new results\,: a proof of
the Burkholder conjecture in probability theory \cite{dede}, a
shorter proof of the Boyadzhiev-de~Laubenfels theorem (concerning
decomposition for group generators, see \cite{haase}), a
characterization for generators of cosine functions
\cite{arendt,cr,haase}, and a characterization of similarities of
$\omega$-accretive operators \cite{lem}.

Let us mention some other consequences.
Assuming $W(A)\subset \overline\Omega$, the inequality
\begin{eqnarray*}
\Vert R(A)\Vert\leq C_{cb}(\Omega)\sup_{z\in \Omega }\|R(z)\|,
\end{eqnarray*}
which holds for all rational functions with matrix values, means that the homomorphism $u_A$ from
the algebra of rational functions bounded on $\Omega $ into the
$C^*$-algebra $\Lh$, defined by $u_A(r) = r(A)$, is completely
bounded with $\Vert u_A\Vert_{cb}\leq
C_{cb}(\Omega )$. A direct application of Paulsen's Theorem (see
\cite{paul} or \cite[Theorem 9.1]{book}) gives
\smallskip

{\it There exists an invertible operator $S\in\Lh$, with
$\Vert S\Vert\,\Vert S^{-1}\Vert\leq
C_{cb}(\Omega )$, such that the domain $\Omega$ is a
complete spectral set for $S\,A\,S^{-1}$.}
\smallskip

We deduce then from a result due to Arveson (see \cite{arveson} or
\cite[Corollary 7.7]{book}) that there exists a
larger Hilbert space $K$ containing $H$ as a subspace
(with the same inner product) and a
normal operator $N$ acting on $K$, with spectrum $\sigma(N)\subset
\partial \Omega $, such that, for all rational functions $r$ bounded in
$\Omega$, we have
\begin{eqnarray*}
r(A)=S^{-1}\,P_H\,r(N)|_H\,S.
\end{eqnarray*}
Here $P_H$ denotes the orthogonal projection from $K$ onto $H$. In
other words, if $W(A)\subset \overline\Omega$,
then $A$ is similar to an operator having a \emph{normal} $\partial
\Omega-$\emph{dilation}. We would like to stress here that our methods give
sharp estimates for the similarity constant
$\Vert S\Vert\,\Vert S^{-1}\Vert$. In particular, we obtain
a similarity constant which is independent of $\Omega$.

\smallskip

Another motivation for our study is that estimates for $C(\Omega)$
and $C_{cb}(\Omega)$ have an interesting potential of applications
in numerical analysis. In computational linear algebra for
instance, the popular Krylov type methods for solving large linear
systems $Ax=b$ are based on polynomial approximations of $A^{-1}$.
We refer to \cite{bb}, where the authors are using some results of
the present paper to improve known error estimates for the GMRES
method. Also, time discretizations of parabolic type P.D.E. use
rational approximations of the exponential. While boundedness is
often sufficient for theoretical considerations, sufficiently good
estimates of our constants are important for numerical
applications.
\smallskip

\subsection{What is this paper about}
The goal of this paper is to present different approaches which
can be used for estimating $C(\Omega)$ and $C_{cb}(\Omega)$. The
obtained estimates are in terms of geometric quantities of
$\Omega$, namely~:

\begin{itemize}
  \item the total variation of $\log (|\sigma - \omega|)$ for a
  bounded convex $\Omega$ and an arbitrary interior point $\omega
  \in \Omega$~;
  \item the angle of a sector contained in the (unbounded) convex
  domain $\Omega$~;
  \item norm estimates for the inverse of the operator used to
  solve the Carl Neumann problem for the double layer potential.
\end{itemize}
We compare these estimates for specific domains several times
in the paper (cf. Remarks~2.1, 3.1, 4.1, 4.3, 4.4 and 4.5). 
The estimates are far better than the general bound $57$, 
except when the domain $\Omega$ is very flat.

The outline of the paper is as follows. The first sections are
based on appropriate integral representations of $r(A)$ or
$R(A)$\,; the positivity (for convex domains) of the double layer
potential plays an important role. In Section~2 we show that
\begin{eqnarray}\label{convb}
C_{cb}(\Omega )\leq 2+\pi +\inf_{\omega  \in\, \Omega }{\rm
TV}(\log|\sigma \!-\!\omega |),
\end{eqnarray}
for every bounded convex domain $\Omega$. Here TV$(\log|\sigma\!-\!\omega |)$ denotes the total variation of
$\log(|\sigma -\omega |)$ as $\sigma $ runs around $\partial \Omega$.
In the unbounded case we obtain in Section~3 the inequality
\begin{eqnarray}\label{convo}
C_{cb}(\Omega )\leq 1+\frac{2}{\pi }\int_\alpha^{\pi /2}\frac{\pi
-x+\sin x}{\sin x}
\, dx,
\end{eqnarray}
if $\Omega $ contains a sector of positive angle $2\alpha $,
$0<\alpha \leq \tfrac{\pi }{2}$. Another representation, based on
the solution of the C.~Neumann problem for the double layer
potential, is given in Section~4. Connections with dilation
theorems are indicated. Section~5 is devoted to the similarity
approach\,; it gives a complete answer for the disk case
$C(D)=C_{cb}(D)=2$ and it is used to show that $\sup_\Omega
C_{cb}(\Omega ,2)=2$, where $C_{cb}(\Omega ,2)$ is defined
similarly as $C_{cb}(\Omega )$ but the supremum is taken now only
over the $2\times2$ matrices. We mention that the methods used in
this last section are genuinely geometric.


\section{The case of a bounded convex domain}\label{sect:2}
Let $\Omega$ be a convex domain of the complex
plane (i.e. $\Omega$ is a convex and open subset of $\C$, $\Omega\neq\emptyset$, $\Omega\neq\C$). On the counterclockwise oriented boundary
$\partial \Omega $, we consider the generic point $\sigma$ of arclength $s$\,; we denote
by $\nu=\tfrac{1}{i}\frac{d\sigma }{ds}$ the unit outward normal (which exists a.e.). Let $A\in\Lh$ be an operator. We introduce the function $\mu (\sigma ,z)$, the half-plane $\Pi _\sigma \supset \Omega $
and the self-adoint operator $\mu(\sigma,A)$ defined by
\begin{eqnarray*}
\mu (\sigma ,z)&=&\frac{1}{\pi }\,\tfrac{\ds d}{\ds ds}\big(\arg(\sigma -z)\big)=
\frac{1}{2\pi}\,(\frac{\nu}{\sigma -z}+\frac{\bar\nu}
{\bar\sigma -\bar z}),\\[2pt]
\Pi _\sigma &=&\{z\,; \mu (\sigma ,z)>0\}=\{z\,; {\Re\,}\bar\nu(\sigma\!-\!z)>0\},\\[2pt]
\mu (\sigma ,A)&=&\tfrac{\ds1}{\ds 2\pi}\,\big(\nu(\sigma\! -\!A)^{-1}+\bar\nu(\bar\sigma -A^*)^{-1}\big)
\end{eqnarray*}
(if $\sigma$ belongs to the resolvent of $A$).
\begin{llemma}
\label{positif}
We assume that the convex domain $\Omega$ contains the spectrum of $A\in\Lh$.
Then the condition $W(A)\subset\Omega\;$ is equivalent to
$$\quad \mu(\sigma,A)>0,\ \forall \sigma\in\partial\Omega.$$
When this condition is satisfied, and if $\Omega$ is bounded and $g$
is a continuous function bounded by 1 on $\partial\Omega$, then we have
\begin{eqnarray*}
\big\|\int_{\partial\Omega} g(\sigma)\,\mu(\sigma,A)\,ds\big\|\leq 2.
\end{eqnarray*}
\end{llemma}
\begin{proof}
We have (setting $w=(\sigma\!-\!A)^{-1}v$)
\begin{eqnarray*}
\mu(\sigma,A)>0&\Longleftrightarrow& \forall \ 0\neq v\in H,\quad \Re\,\bar\nu\langle
(\bar\sigma\!-\!A^*)^{-1}v, v\rangle>0\\
&\Longleftrightarrow& \forall \ 0\neq w\in H,\quad \Re\,\bar\nu\langle
(\sigma\!-\!A)w, w\rangle>0\\
&\Longleftrightarrow& W(A)\subset \Pi_{\sigma}.
\end{eqnarray*}
The equivalence follows from the convexity property $\Omega= \cap_{\sigma\in\partial\Omega}\Pi_{\sigma}$.\smallskip

We deduce from the Cauchy formula that $\int_{\partial\Omega}\mu(\sigma,A)\,ds=2$.
We set $\Gamma=\int_{\partial\Omega}g(\sigma)\,\mu(\sigma,A)\,ds$ and consider
two norm one vectors $u,v\in H$. From the positivity of $\mu(\sigma,A)$ we have
\begin{eqnarray*}
|\langle\mu(\sigma,A)\,u,v\rangle|\leq\tfrac{1}{2}\,\langle\mu(\sigma,A)\,u,u\rangle+
\tfrac{1}{2}\,\langle\mu(\sigma,A)\,v,v\rangle;
\end{eqnarray*}
therefore
\begin{eqnarray*}
|\langle \Gamma\,u,v\rangle|&\leq &
\tfrac{1}{2}\,\int_{\partial\Omega}\langle\mu(\sigma,A)\,u,u\rangle\,ds+
\tfrac{1}{2}\,\int_{\partial\Omega}\langle\mu(\sigma,A)\,v,v\rangle\,ds\\
 & = &
\|u\|^2+\|v\|^2=2.
\end{eqnarray*}
The result follows from  $\|\Gamma\|=\sup\{|\langle
\Gamma\,u,v\rangle|\,;\|u\|=\|v\|=1\}$.
\end{proof}
\noindent{\it Remark.} With an easy modification this proof is also valid if $g$ is a matrix valued function.
\smallskip

For the remaining part of this section we assume that the domain $\Omega$ is bounded and that the origin
$0$ belongs to $\Omega$.
Now we introduce the angles $\theta$ and $\varphi$ as described below.
\begin{center}
\centerline{\psfig{figure=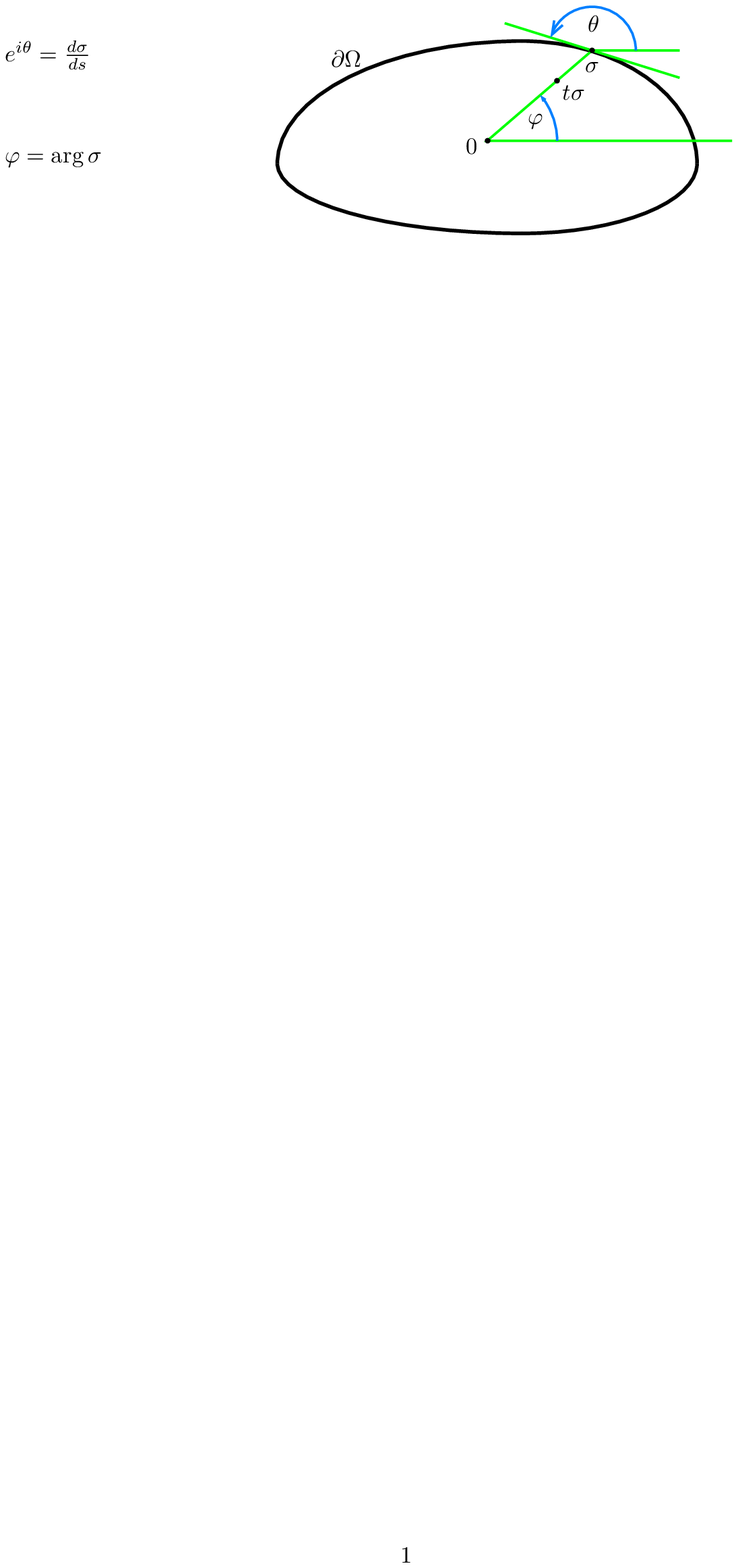,height=3.7cm}}
\end{center}
\medskip

\noindent We can choose $\theta$ and $\varphi$ such that
$0<\theta-\varphi<\pi$. In order to avoid some technical
difficulties, we initially assume that
\begin{eqnarray} \label{regular}
\theta\ \hbox{ is a $C^1$ function of $s$ and } \theta'(s)>0, \forall
\,s\in[0,L].
\end{eqnarray}
Here $L$ is the length of $\partial\Omega$.
This assumption allows to consider also $s$, $\sigma$ and $\varphi$  as $C^1$ functions of $\theta$.

\begin{llemma}
\label{thmu} Let $r$ be a rational function bounded in a domain $\Omega $ satisfying (\ref{regular}). Then we have
\begin{eqnarray}\label{representation}
r(z)=\int_{\partial \Omega } r(\sigma )\mu (\sigma ,z)\,ds +
\int_0^{2\pi } Jr(\sigma ,\bar z) \,d\theta ,\qquad\forall z\in\Omega,
\end{eqnarray}
where $\ds\quad Jr(\sigma ,\bar z):=
\frac{1}{\pi }\,\int_0^1\frac{e^{2i\theta}\sigma(\bar\sigma\!-\!\bar z)\,r(t\sigma)}
{\big((t\!-\!1)\sigma+e^{2i\theta}(\bar\sigma\!-\!\bar z)\big)^2}\, dt$.
\smallskip

\noindent Furthermore, if $|r|\leq1$ in $\Omega$, then we have the estimate
\begin{eqnarray}\label{supiteta}
\vert Jr(\sigma,\bar z)\vert \leq \tfrac{1}{2}\!+\!|\cot(\theta\!-\!\varphi )|,\qquad \forall \,z\in \Pi_{\sigma}.
\end{eqnarray}
\end{llemma}
\begin{proof}
a) For a given point $z\in\Omega$ we introduce ($\sigma$ and $\theta$ being functions of $s$)
\begin{eqnarray*}
u(t,s)&:=& (t\!-\!1)\sigma+e^{2i\theta}(\bar\sigma\!-\!\bar z),\\
v(t,s)&:=& \frac{\sigma\,r(t\sigma)}{u},\quad w(t,s):=\frac{t\,r(t\sigma)}{u}\,\tfrac{\ds d\sigma}{\ds ds}.
\end{eqnarray*}
We first notice that
\begin{eqnarray*}
\frac{\partial}{\partial s}(u\,v)-\frac{\partial}{\partial t}(u\,w)=0,
\qquad\qquad d\sigma=e^{2i\theta}d\bar\sigma,\\
\frac{\partial u}{\partial s}=t\frac{d\sigma}{ds}+2\,i\,e^{2i\theta}(\bar\sigma\!-\!\bar z)
\frac{d\theta}{ds},\qquad\frac{\partial u}{\partial t}=\sigma.\quad
\end{eqnarray*}
Therefore
\begin{eqnarray*}
u\,\big(\frac{\partial w}{\partial t}-\frac{\partial v}{\partial s}\big)=
v\,\frac{\partial u}{\partial s}-w\,\frac{\partial u}{\partial t}=
2\,i\,e^{2i\theta}\,v\,(\bar \sigma\!-\!\bar z)\frac{d\theta}{ds},\\
\hbox{and}\qquad \frac{\partial w}{\partial t}-\frac{\partial v}{\partial s}=2\,i\,
\frac{e^{2i\theta}\sigma(\bar\sigma\!-\!\bar z)\,r(t\sigma)}
{\big((t\!-\!1)\sigma+e^{2i\theta}(\bar\sigma\!-\!\bar z)\big)^2}\,\frac{d\theta}{ds}.\hskip2cm
\end{eqnarray*}
We also have, if $t\in[0,1]$,
\begin{eqnarray*}
\Im(e^{-i\theta}u)&=&(t\!-\!1)\Re (-\bar\nu\sigma)+\Re \nu(\bar\sigma\!-\!\bar z)\\&=&(1\!-\!t)\Re \bar\nu(\sigma\!-\! 0)+\Re \bar\nu(\sigma\!-\! z)>0;
\end{eqnarray*}
indeed $0$ and $z$ are in $\Omega\subset\Pi_{\sigma}$. Thus the function $u(t,s)$ does not vanish for $(t,s)\in[0,1]\times[0,L]$\,; this shows that the functions $v$, $\frac{\partial v}{\partial s}$, $w$,
$\frac{\partial w}{\partial t}$ are well defined and continuous on this set. Therefore we can write
\begin{eqnarray*}
2\,\pi\,i\,\int_{0}^{2\pi}Jr(\sigma,\bar z)\,d\theta&=&\int_{0}^{L}\!\!\int_{0}^1\big(\frac{\partial w}{\partial t}-\frac{\partial v}{\partial s}\big)\,dt\,ds\\
&=&\int_{0}^{L}\!\!\int_{0}^1\frac{\partial w}{\partial t}\,dt\,ds-
\int_{0}^{1}\!\!\int_{0}^L\frac{\partial v}{\partial s}\,ds\,dt\\
&=&\int_{0}^{L} w(1,s)\,ds=\int_{\partial\Omega} \frac{r(\sigma)}{\bar \sigma-\bar z}\,d\bar \sigma.
\end{eqnarray*}
We deduce that
\begin{eqnarray*}
\int_{\partial \Omega } r(\sigma )\mu (\sigma ,z)\,ds +
\int_0^{2\pi } Jr(\sigma ,\bar z) \,d\theta=\frac{1}{2\pi i}\,\int_{\partial\Omega} \frac{r(\sigma)}{\sigma-z}\,d\sigma,
\end{eqnarray*}
and then (\ref{representation}) follows from the Cauchy formula.
\medskip

b) We remark that $Jr(\sigma,\bar z)$ is antiholomorphic for $z\in\Pi_{\sigma}$, bounded and continuous in $\bar\Pi_{\sigma}$ and vanishing as $z\to \infty$. Using the maximum principle we obtain for $z\in\Pi_{\sigma}$
\begin{eqnarray*}
|Jr(\sigma,\bar z)|&\leq&\sup_{\zeta\in\partial\Pi_{\sigma}} |Jr(\sigma,\bar \zeta)|
=\sup_{x\in\R} |Jr(\sigma,\bar \sigma+xe^{-i\theta})|\\
&\leq&\sup_{x\in\R}\int_{0}^1\frac{|\sigma|\,|x|}{\pi\,|(t-1)\sigma-xe^{i\theta}|^2}\,dt\\
&\leq&\max_{\varepsilon=\pm1}\int_{0}^\infty \frac{d\tau}{\pi\,|\tau e^{i\varphi}-\varepsilon e^{i\theta}|^2}
 =\frac{\max(\theta -\varphi,\pi -\theta+\varphi)}{\pi\,\sin(\theta-\varphi)}.
\end{eqnarray*}
We have used above the change of variables $\tau=(t\!-\!1)|\sigma|/|x|$.
Now we remark that the inequality $\max(u,\pi\!-\!u)\leq\frac{\pi}{2}\sin u+\pi|\cos u|$ holds for any $u=\theta\!-\!\varphi\in(0,\pi)$, which shows that
\begin{eqnarray*}
|Jr(\sigma,\bar z)|\leq\frac{1}{2}+|\cot (\theta\!-\!\varphi)|.
\end{eqnarray*}
\end{proof}

\begin{theorem}
\label{noflt}
For all bounded convex domains $\Omega$ we have the estimate
\begin{eqnarray*}
C(\Omega )\leq C_{cb}(\Omega )\leq 2+\pi +\phi_{_\Omega},\quad\hbox{ where }
\phi_{_\Omega}:=\inf_{\omega\in\Omega} TV(\log|\sigma\!-\!\omega|).
\end{eqnarray*}
\end{theorem}
\begin{proof}
Without loss of generality we can assume that $\omega=0$ is the origin of the complex plane.
We assume initially that (\ref{regular}) is satisfied.
Let $A$ be an operator with $\overline{W(A)}\subset \Omega$ and $r$ a rational function bounded by 1 in $\Omega$.
We can replace $z$ by $A$ and $\bar z$ by $A^*$ in the proof of the previous lemma. We obtain
\begin{eqnarray*}
2\,\pi\,i\,\int_{0}^{2\pi}Jr(\sigma,A^*)\,d\theta=
\int_{\partial\Omega} r(\sigma)\,(\bar \sigma\!-\!A^*)^{-1}d\bar \sigma,
\end{eqnarray*}
and then, from the Cauchy formula,
\begin{eqnarray*}
r(A)=\int_{\partial \Omega } r(\sigma )\mu (\sigma ,A)\,ds +
\int_0^{2\pi } Jr(\sigma ,A^*) \,d\theta.
\end{eqnarray*}
From Lemma \ref{positif} we know that $\|\int_{\partial \Omega } r(\sigma )\mu (\sigma ,A)\,ds\|\leq2$.
Using the von Neumann inequality for the half-plane $\Pi_{\sigma}$ we have
\begin{eqnarray*}
\| Jr(\sigma ,A^*) \|\leq\sup_{z\in \Pi_{\sigma}} |Jr(\sigma,\bar z)|\leq  \tfrac{\ds1}{\ds2}+|\cot (\theta\!-\!\varphi)|.
\end{eqnarray*}
Therefore
\begin{eqnarray*}
\|r(A)\|\leq 2+\pi+ \int_0^{2\pi } |\cot (\theta\!-\!\varphi)| \,d\theta=2+\pi+ \int_0^{2\pi } |\cot (\theta\!-\!\varphi)| \,d\varphi;
\end{eqnarray*}
indeed $\int_0^{2\pi } |\cot (\theta\!-\!\varphi)|
\,d(\theta-\varphi)=0$. Setting $\rho(\varphi):=|\sigma|$, we have
$\cot(\theta\!-\!\varphi)=\rho'(\varphi)/\rho(\varphi)$ and thus
\begin{eqnarray*}
\|r(A)\|\leq 2+\pi+ \int_0^{2\pi } \Big|\frac{\rho'(\varphi)}{\rho(\varphi)}\Big|\,d\varphi=
2+\pi+ TV(\log \rho).
\end{eqnarray*}
We obtain the estimate $C(\Omega )\leq 2+\pi +\phi_{_\Omega}$.\smallskip

Now we remark that the von Neumann inequality,
as well as Lemma \ref{positif}, are also valid for rational functions with matrix values\,; therefore we obtain in the same way $C_{cb}(\Omega )\leq 2+\pi +\phi_{_\Omega}$. The use of Lemma \ref{reg} hereafter (with $K=\overline{W(A)}$) allows to extend the results to all bounded domains $\Omega$.
\end{proof}

\noindent{\bf Remark 2.1.} For instance, if $\Omega $ is an ellipse
with major
axis $2a$ and minor axis $2b$, then
we have $\phi_{_\Omega} =4\log\frac{a}{b}$ \,;
$\phi_{_\Omega} =0$ if $\Omega $ is a disk. For some domains $\Omega$ it is
dificult to compute exactly $\phi_{_\Omega}$. For these domains,
the following geometrical estimate given in \cite{crac} can be useful
\begin{eqnarray*}
\phi_{_\Omega}\leq2\,\log\frac{\tau_{_\Omega}^4}{\tau_{_\Omega}^2-2},\qquad\hbox{where\quad}
\tau_{_\Omega}:=\min_{\omega\in\Omega}\frac{\max\{|\sigma-\omega|\,;\sigma\in\partial\Omega\}}
{\min\{|\sigma-\omega|\,;\sigma\in\partial\Omega\}}.
\end{eqnarray*}
The quantity $\tau_{_{\Omega}}$ can be considered as a rate of
flatness of $\Omega$.

We obtain $C_{cb}(\Omega) \leq 14.4$ if $\tau_\Omega \leq 10$ and 
$C_{cb}(\Omega) \leq 40$ if $\tau_\Omega \leq 6000$. 
The obtained estimate for $C_{cb}(\Omega)$ is 
better than the general
bound $57$ given in \cite{crac}
for $\tau_{_{\Omega}}\leq 400\,000$. 
\begin{llemma}\label{reg} We assume that $K$ is a compact subset of a  bounded convex
domain $\Omega$ and that $0\in\Omega$.
Then, for all $\varepsilon>0$, there exists a bounded convex
domain $\Omega'$ satisfying (\ref{regular}) such that
\begin{eqnarray*}
K\subset\Omega'\subset\Omega,\quad\hbox{and\quad}TV_{\partial\Omega'}(\log|\sigma|)\leq
TV_{\partial\Omega}(\log|\sigma|)+\varepsilon.
\end{eqnarray*}
\end{llemma}
\begin{proof} In this lemma we use the notation $r(\varphi)=|\sigma|^{-1}=1/\rho(\varphi)$,
where $\varphi$ is the argument of the boundary point $\sigma$. The convexity of $\Omega$ is
equivalent to
\begin{quote}
 $r$ is a continuous $2\pi$-periodic function on $\R$, and $r+r''$ is a
 positive measure.
\end{quote}

\noindent By a standard mollifier technique we can
find a sequence of $2\pi$-periodic functions $r_{n}\in C^\infty(\R)$
such that
\begin{quote}
$r_{n}\to r$  uniformly in $\R$, \quad $\min_{x}(r_{n}(x)+r_{n}''(x))\geq 0$ and \quad$r'_{n}\to r'$
in $L^1(0,2\pi)$.
\end{quote}

\noindent By adding (if needed) to $r_{n}$ a positive constant,
we can assume
that $r_{n}\geq r$ and $\min_{x}(r_{n}(x)+r_{n}''(x))>0$. We set
\begin{eqnarray*}
\Omega_{n}:=\{\rho\,e^{i\varphi}\,; 0\leq\rho<1/r_{n}(\varphi), \varphi\in[0,2\pi]\}.
\end{eqnarray*}
It is clear that $\Omega_{n}$ is open, strictly convex, satisfies (\ref{regular}), that
$\Omega_{n}\subset\Omega$ and that, for $n$ large enough, $K\subset \Omega_{n}$.
Furthermore, we have
\begin{eqnarray*}
TV_{\partial\Omega_{n}}(\log|\sigma|)=\int_{0}^{2\pi}\frac{|r_{n}'(\varphi)|}{r_{n}(\varphi)}\,
d\varphi\to\int_{0}^{2\pi}\frac{|r'(\varphi)|}{r(\varphi)}\,d\varphi=TV_{\partial\Omega}(\log|\sigma|),
\end{eqnarray*}
as $n$ tends to infinity.
The Lemma follows by choosing $\Omega'=\Omega_{n}$ with a sufficiently large value of $n$.
\end{proof}
\section{Unbounded domains}\label{sect:3}

We turn now to the case of an unbounded convex domain $\Omega\neq \C$.
We assume in this section that $\Omega$ contains a convex sector of angle $2\alpha>0$.
Since the constants are only dependent on the shape of $\Omega$, we can assume, without loss of generality, that the value $\alpha$ is maximal and that
\begin{eqnarray*}
S_{\alpha}:=\{z\in\C\,; z\neq0, |\arg z|<\alpha\}\subset\Omega\subset\{z\in\C\,; \Re z>0\}.
\end{eqnarray*}

We use the same notations as in  the previous section.
Let $A\in\Lh$ be an operator with $\overline{W(A)}\subset\Omega$. We remark that
\begin{eqnarray*}
\int_{\partial\Omega}\mu(\sigma,A)\,ds=2-2\,\alpha/\pi.
\end{eqnarray*}
Indeed, if we choose $R>\|A\|$ and set $\Omega_{R}=\{z\in\Omega\,;|z|<R\}$,
$C_{R}=\partial\Omega_{R}\setminus(\partial\Omega\cap\partial\Omega_{R})$,
then we have
\begin{eqnarray*}
\int_{\partial\Omega_{R}}\mu(\sigma,A)\,ds=2
\end{eqnarray*}
and
\begin{eqnarray*}
\int_{C_{R}}\mu(\sigma,A)\,ds=\int_{C_{R}}\frac{1}{\pi\,R}\,ds+O(\frac{\|A\|}{R})=
\frac{2\,\alpha}{\pi}+O(\frac{1}{R}).
\end{eqnarray*}
Thus
\begin{eqnarray*}\label{muunb}
\int_{\partial\Omega}\mu(\sigma,A)\,ds=\lim_{R\to\infty}\int_{\partial\Omega\cap\partial\Omega_{R}}\mu(\sigma,A)\,ds=2-2\,\alpha/\pi.
\end{eqnarray*}
This implies that Lemma \ref{positif} is still valid for unbounded domains,
but now with a better estimate~: if $g$ is a continuous
function bounded by $1$ on $\partial\Omega$, then
\begin{eqnarray*}
\Big\|\int_{\partial\Omega}g(\sigma)\,\mu(\sigma,A)\,ds\Big\|\leq2-2\,\alpha/\pi.
\end{eqnarray*}

We have now to modify Lemma \ref{thmu} as follows.
\begin{llemma}
\label{lemu} We assume that $\theta\in C^1(\R)$ satisfies $\theta'(s)>0$
for all real $s$. Let $r$ be a rational function bounded in $\Omega $.
Then we have
\begin{eqnarray}\label{decos}
r(z)=\int_{\partial \Omega } r(\sigma )\mu (\sigma ,z)\,ds +
\int_{\pi+\alpha}^{2\pi-\alpha} Kr(\sigma ,\bar z) \,d\theta ,\qquad\forall z\in\Omega,
\end{eqnarray}
where $\ds\quad Kr(\sigma ,\bar z):=-\frac{1}{\pi }\,\int_{0}^{\infty}
  r(\sigma\!+\!t)\, \frac{ e^{2i\theta }(\bar \sigma -\bar z )}{(t
+e^{2i\theta } (\bar \sigma \!-\! \bar z ))^2}\ dt.$
\smallskip

\noindent Furthermore, if $|r|\leq1$ in $\Omega$, then we have the estimate
\begin{eqnarray}\label{supitete}
|Kr(\sigma ,\bar z)|\leq \frac{\max(\theta -\pi,2\pi -\theta)}
{\pi \,\sin\theta} ,\qquad \forall \,z\in \Pi_{\sigma}.
\end{eqnarray}
\end{llemma}
\begin{proof} a) We first remark that if $r=1$ is constant, then
$Kr(\sigma,\bar z)=-1/\pi$ and (\ref{decos}) is satisfied.
It is thus sufficient to suppose
$r(\infty)=0$ in the proof of (\ref{decos}).

We note that the angle $\theta $ of the tangent vector in the boundary point $\sigma \in\partial\Omega$ is now running in the interval $(\pi +\alpha ,2\pi -\alpha )$. For a fixed $z\in\Omega$ we set
\begin{eqnarray*}
u(t,s):= t+e^{2i\theta}(\bar\sigma\!-\!\bar z),\quad
v(t,s):= \frac{r(\sigma\!+\!t)}{u},\quad w(t,s):=\frac{r(\sigma\!+\!t)}{u}\,\tfrac{\ds d\sigma}{\ds ds}.
\end{eqnarray*}
We remark that $\Im (e^{-i\theta}u)=-t\,\sin\theta+\Re\bar\nu(\sigma-z)>0$. Thus
the function $u$ does not vanish in $[0,\infty[\times\R$ and there exists a constant $c=c(z)>0$ such that on this set $|u(t,s)|\geq c\,(1\!+\!t\!+\!|s|)$. We have
\begin{eqnarray*}
u^2\,\big(\frac{\partial w}{\partial t}-\frac{\partial v}{\partial s}\big)&=&
r'(\sigma\!+\!t)\,\frac{d\sigma}{ ds}\,u
-r(\sigma\!+\!t)\,\frac{ d\sigma}{  ds}\,\tfrac{\ds\partial u}{\ds\partial
  t}\\
&-& 
r'(\sigma\!+\!t)\,\frac{  d\sigma}{ ds}\,u
+r(\sigma\!+\!t)\,\frac{ \partial u}{\ \partial s}\\
&=&2\,i\,r(\sigma\!+\!t)\,e^{2i\theta}(\bar\sigma\!-\!\bar z)\,\frac{d\theta}{ds}.
\end{eqnarray*}
Thus
\begin{eqnarray*}
2\,\pi\,i\,Kr(\sigma ,\bar z)\,\frac{d\theta}{ds}&=&\int_{0}^\infty\frac{\partial v}{\partial s}(t,s)\,dt
-\int_{0}^{\infty}\frac{\partial w}{\partial t}\,dt\\
 &=&\frac{r(\sigma)}{e^{2i\theta}(\bar\sigma-\bar z)}\,
\frac{ d\sigma}{ds}+\int_{0}^\infty\frac{\partial v}{\partial s}(t,s)\,dt.
\end{eqnarray*}
Noticing that $|\frac{\partial v}{\partial s}(t,s)|\leq
C\,(1\!+\!t\!+\!|s|)^{-3}$, which justifies the use of Fubini's theorem,
we have
\begin{eqnarray*}
2\,\pi\,i\,\int_{\pi+\alpha}^{2\pi-\alpha}Kr(\sigma ,\bar z)\,d\theta&=&\int_{\partial\Omega} \frac{r(\sigma)}{\bar\sigma-\bar z}\,d\bar\sigma+
\int_{-\infty}^{+\infty} \!\!\int_{0}^\infty\frac{\partial v}{\partial s}(t,s)\,dt\,ds\\
&=&\int_{\partial\Omega} \frac{r(\sigma)}{\bar\sigma-\bar z}\,d\bar\sigma+
\int_{0}^\infty\!\!\int_{-\infty}^{+\infty}\frac{\partial v}{\partial
  s}(t,s)\,dt\,ds\\ &=&
\int_{\partial\Omega} \frac{r(\sigma)}{\bar\sigma-\bar z}\,d\bar\sigma.
\end{eqnarray*}
The relation (\ref{decos}) follows now from the Cauchy formula.
\medskip

b) Using the maximum principle we obtain, for $z\in\Pi_{\sigma}$,
\begin{eqnarray*}
|Kr(\sigma,\bar z)|&\leq&\sup_{\zeta\in\partial\Pi_{\sigma}} |Kr(\sigma,\bar \zeta)|
=\sup_{x\in\R} |Kr(\sigma,\bar \sigma+xe^{-i\theta})|\\
&\leq&\sup_{x\in\R}\int_{0}^\infty\frac{|x|}{\pi\,|t-xe^{i\theta}|^2}\,dt\\
&\leq&\max_{\varepsilon=\pm1}\int_{0}^\infty \frac{d\tau}{\pi\,|\tau-\varepsilon e^{i\theta}|^2}
 =\frac{\max(\theta -\pi,2\pi -\theta)}{\pi\,\sin\theta}.
\end{eqnarray*}
\end{proof}

\begin{theorem}
\label{unbdeux}
We assume that the convex domain $\Omega \neq \C$ contains a sector of
angle $2\alpha $, with $0<\alpha \leq \pi /2$. Then we have the estimate
\begin{eqnarray*}
C(\Omega)\leq C_{cb}(\Omega) \leq  1+\frac{2}{\pi}\int_\alpha ^{\pi /2} \frac{\pi -x+\sin x }{\sin x}\,dx.
\end{eqnarray*}
\end{theorem}
\begin{proof} As for Theorem \ref{noflt}, it is sufficient to look at 
the estimate of $C(\Omega)$
and we can  assume that $\theta\in C^1(\R)$ satisfies $\theta'(s)>0$ for all
real $s$. Let $r$ be a rational function bounded by 1 in $\Omega $. We
deduce from (\ref{decos}) that
\begin{eqnarray*}
r(A)=\int_{\partial \Omega } r(\sigma )\mu (\sigma ,A)\,ds +
\int_{\pi+\alpha}^{2\pi-\alpha} Kr(\sigma ,A^*) \,d\theta.
\end{eqnarray*}
We have seen that
\begin{eqnarray*}
\Vert \int_{\partial \Omega } r(\sigma )\mu (\sigma ,A)\,ds \Vert \leq
2-2\,\alpha/\pi
\end{eqnarray*}
and
\begin{eqnarray*}
\Vert  Kr(\sigma ,A^*) \Vert \leq \frac{\max(\theta -\pi,2\pi -\theta)}{\pi \,\sin\theta}.
\end{eqnarray*}
Alltogether that gives
\begin{eqnarray*}
\Vert r(A)\Vert \leq 1+\frac{2}{\pi }\int_\alpha^{\pi /2} \frac{\pi -x+\sin x}{\sin x}\,dx
\end{eqnarray*}
and the theorem follows.
\end{proof}

\noindent{\bf Remark 3.1.} We obtain 
$C_{cb}(\Omega) \leq 6$ if $\alpha \leq \pi/20$ ; 
$C_{cb}(\Omega) \leq 40$ if $\alpha \leq 10^{-8}$ and
$C_{cb}(\Omega) \leq 57$ if $\alpha \leq 1.3\cdot 10^{-12}$.

\section{A potential-theoretic approach}\label{sect:4}
\subsection{The Carl Neumann problem}\label{subsect:41}

Let $\Omega \neq \emptyset$ be a bounded convex domain of the complex
plane.
Given a function $r$, continuous on $\overline\Omega $ and
harmonic in $\Omega
$, the C. Neumann problem \cite{cneu} for the double layer potential on
$\partial
\Omega
$ is the following : find a function $g\in C(\partial \Omega )$ such
that
\begin{eqnarray}\label{cne}
\forall \,z\in \Omega,\quad r(z)=\frac{1}{2}\int_{\partial \Omega
}g(\sigma
)\,\mu (\sigma ,z)\,ds.
\end{eqnarray}
Taking the limit as the point $z$ tends to the
boundary $\partial \Omega$ ,  the problem (\ref{cne}) is
equivalent to
\begin{eqnarray*}
\forall \,z\in \partial \Omega ,\quad
r(z)=\frac{1}{2}\,\big(g(z)+\int_{\partial \Omega_{z} }g(\sigma )\,\mu
(\sigma ,z)\,ds\big) ,
\end{eqnarray*}
where
\begin{eqnarray*}
\partial \Omega_{z}=\partial \Omega\setminus\{z\}.
\end{eqnarray*}
This relation can be written (by considering restrictions to $\partial \Omega $)
\begin{eqnarray}\label{neub}
r=\tfrac{1}{2}(I+P)g,\quad \hbox{where}\quad Pg(z)=\frac{1}{\pi
}\int_{\partial \Omega_{z} }g(\sigma )\,d\arg(\sigma -z).
\end{eqnarray}
Clearly $P$ is a linear operator acting from
$C(\partial \Omega )$ into itself.
A harmonic function is uniquely defined by its restriction on
the boundary $\partial\Omega $, and any continuous function on this
boundary is the trace of such a function. Therefore the invertibility of
the operator $I+P$ in $C(\partial \Omega )$ implies existence and
uniqueness for (\ref{cne}). It is known that this operator is
effectively
invertible, since we have
assumed $\Omega $ bounded and convex (see for instance the monograph
\cite{kral}).
\smallskip

We will mainly restrict our attention to rational
functions, bounded in $\Omega $, and introduce the constant
\begin{eqnarray*}
C_N(\Omega )=\sup\{2\,\Vert (I\!+\!P)^{-1}r\Vert _{L^\infty(\partial
\Omega )} \,; r\
\hbox{rational function},\ |r|\leq1\ \hbox{ in  } \Omega \}.
\end{eqnarray*}
In this definition  we assume that $r$ acts from $\C$ into $\C$,
but the constant would be unchanged by considering matrix-valued
rational functions.\smallskip

Thus, if $R$ is a matrix-valued rational function satisfying $\|R\|\leq1$ in $\Omega$, we  have, setting
$G=2(I+P)^{-1}R$,
\begin{eqnarray*}
\Vert G\Vert _{L^\infty(\partial \Omega )}\leq C_N(\Omega )\quad\hbox{and}
\quad
R(z)=\frac{1}{2}\,\int_{\partial \Omega } G(\sigma )\mu (\sigma ,z)\,ds.
\end{eqnarray*}
Comparing the holomorphic and antiholomorphic parts we deduce that,
for some complex constant $c$,
\begin{eqnarray*}
R(z)=\frac{1}{4\,\pi}\,\int_{\partial \Omega } G(\sigma )\,\frac{\nu}{\sigma-z}\,ds-c
\quad\hbox{and} \quad 0=\frac{1}{4\,\pi}\,\int_{\partial \Omega } G(\sigma )\,\frac{\bar\nu}{\bar\sigma-\bar z}\,ds+c.
\end{eqnarray*}
If $A\in
\Lh$ is an operator with $\overline {W(A)}\subset \Omega $, it is licit to
replace $z$ by $A$ and $\bar z$ by $A^*$ in the previous relations. After
adding we obtain that
\begin{eqnarray*}
R(A)=\frac{1}{2}\,\int_{\partial \Omega } G(\sigma )\otimes \mu (\sigma
,A)\,ds,
\end{eqnarray*}
and, using Lemma \ref{positif}, that
\begin{eqnarray*}
\Vert R(A)\Vert \leq \sup_{z\in \Omega } \Vert G(z)\Vert .
\end{eqnarray*}
Therefore we deduce $C_{cb}(\Omega)\leq C_{N}(\Omega)$.
\medskip

In order to estimate the constant $C_{N}(\Omega)$ it is interesting to
introduce another constant $D_N(\Omega)$ defined as
\begin{eqnarray*}
\sup\left\{2\inf_{c\in\C}\Vert(I\!+\!P)^{-1}r-c\Vert_{L^\infty(\partial\Omega)}\,; r
\hbox{ rational function}, |r|\leq1 \hbox{ in  } \Omega \right\}.
\end{eqnarray*}
It is clear that $D_{N}(\Omega)\leq C_{N}(\Omega)$,  $P1=1$ and $\|Pr\|_{L^\infty(\partial\Omega )}
\leq \| r\|_{L^\infty(\partial\Omega )}$. From the relation
\begin{eqnarray*}
2(I\!+\! P)^{-1}r=r+(I\!-\! P)\big((I\!+\! P)^{-1} r-c\big), \quad\forall c\in\C,
\end{eqnarray*}
we deduce $ C_N(\Omega )\leq 1+D_N(\Omega)$. Therefore we have
\begin{eqnarray}
C_{cb}(\Omega)\leq C_N(\Omega )\leq 1+D_N(\Omega).
\end{eqnarray}

\medskip

For an operator $M\in {\cal L}(C(\partial \Omega ))$ we introduce the
norm and semi-norm
\begin{eqnarray*}
\Vert M\Vert _\infty& =&\sup\, \{\Vert M f\Vert _{L^\infty (\partial \Omega)}\,;
f\in C(\partial \Omega ),\ |f|\leq 1 \hbox{  in  }\Omega\} ,\\
\Vert M\Vert _{\rm osc} &=&\sup\, \{\inf_{c\in\C}\Vert M f\!-\!c\Vert _{L^\infty (\partial \Omega)}\,;
f\in C(\partial \Omega ),\ |f|\leq 1 \hbox{  in  }\Omega\} .
\end{eqnarray*}
Recall that osc\,$(f):=2\,\inf_{c\in \C}\Vert f-c\Vert  _{L^\infty
(\partial \Omega)}$ is called the \emph{oscillation} of $f$ on $\partial
\Omega$.\medskip

We clearly have
$$C_N(\Omega )\leq 2\, \Vert (I+P)^{-1}\Vert_\infty  \quad  ; \quad D_N(\Omega )\leq 2\,\Vert(I+P)^{-1}\Vert_{\rm osc}.$$

\medskip

\noindent {\bf Remark 4.1.} If $\Omega =D$ is a disk of center $\omega
$, then
simple calculations show that $(Pf)(z)=f(\omega )$ for all harmonic
functions $f$ in $\Omega $, with $f\in C(\overline{\Omega})$.
We get $(2(I\!+\! P)^{-1}f)(z)=2f(z)-f(\omega )$ ; therefore
\begin{eqnarray*}
C_{cb}(D)\leq C_N(D)=3\qquad\hbox{and  }\quad D_N(D)=2.
\end{eqnarray*}
The estimate for the first constant is not optimal, since we will see in
Theorem \ref{thdisk} that
$C(D)=C_{cb}(D)=2$.

\medskip
\noindent {\bf Remark 4.2.} It is known \cite{Sch} that, if $\Omega$
is not a triangle nor a quadrilateral, then $\Vert P\Vert _{\rm osc}<1$.
This gives the estimate $D_N(\Omega )\leq 2(1-\Vert P\Vert _{\rm osc})
^{-1}$. Furthermore, using the notations
\begin{eqnarray*}
R_\Omega &:=&\sup\{\hbox{radii of circles which intersect $\partial
\Omega $ in at least 3 points}\},\\
L_\Omega &:=& \hbox{perimeter of }\Omega ,
\end{eqnarray*}
we have $\Vert P\Vert _{\rm osc}\leq 1-\frac{L_\Omega }{2\pi R_\Omega
}$ (some smoothness assumptions are mentioned in \cite{Sch}, but we think
that they can be avoided). We obtain the estimates
\begin{eqnarray*}
D_N(\Omega )\leq \frac{4\pi R_\Omega }{L_\Omega }\qquad\hbox{and  }\quad
C_{cb}(\Omega )\leq 1+\frac{4\pi R_\Omega }{L_\Omega }.
\end{eqnarray*}
These estimates are not useful if $\Omega$ is a polygon since then
$R_\Omega =+\infty $. For an ellipse $(E)$ with major axis $2a$ and
minor axis $2b$, it can be computed that $R_E=a^2/b$\,; therefore
$D_N(E)\leq 2\pi a/b$.

\medskip
\noindent {\bf Remark 4.3.} The modern proof \cite{kral} of the
invertibility of
$I\!+\! P$, which works for any bounded convex domain, follows
from the inequality $\Vert P^2\Vert _{\rm osc}< 1$. Then, from the
relation
$(I\!+\! P)^{-1}=(I\!-\! P)(I\!-\! P^2)^{-1}$, we deduce
\begin{eqnarray*}
D_N(\Omega )\leq 2\,\Vert I\!-\! P\Vert _\infty (1\!-\!
\Vert P^2\Vert _{\rm osc})^{-1}.
\end{eqnarray*}
We are not acquainted with a translation of this inequality in simple
 geometric terms. This was done in another approach developped in
 \cite{dede} which is based on the study of the iterate $P^3$ and which provides
the estimate
 \begin{eqnarray*}
 C_N(\Omega )\leq 3+\Big(\frac{2\pi \delta ^2}{|\Omega
 |}\Big)^3,\quad \delta \hbox{ diameter of }\Omega ,\quad |\Omega |
 \hbox{ area of }\Omega .
 \end{eqnarray*}

\subsection{Extension to unbounded domains}\label{subsect:42}
The previous developments admit an extension to unbounded domains.

\begin{theorem}
\label{unbd}
We assume that the convex domain $\Omega \neq \C$ contains a sector of
angle
$2\alpha $, with $0<\alpha \leq \pi /2$. Then the C. Neumann problem
(\ref{cne}) is well posed and we have the estimates
\begin{eqnarray*}
C_N(\Omega) \leq \frac{\pi}{\alpha }\quad\hbox{and}\quad C_{cb}(\Omega
)\leq \frac{\pi -\alpha }{\alpha  }.
\end{eqnarray*}
\end{theorem}
\begin{proof}

a) As above, we have
\begin{eqnarray*}
\Vert P\Vert _\infty
\leq \frac{1}{\pi }\sup_{z\in \partial \Omega }\int_{\partial \Omega }
d\arg(\sigma -z)=\frac{\pi -2\alpha }{\pi }<1.
\end{eqnarray*}
We deduce that $(I+P)$ is invertible from $C({\partial\Omega} )$ into itself.
Thus the C.~Neumann problem has a unique solution and we have
\begin{eqnarray*}
\Vert r\Vert _{L^\infty (\partial \Omega )}\geq \tfrac{1}{2} (1-\Vert
P\Vert
_\infty) \,\Vert
g\Vert  _{L^\infty (\partial \Omega )},
\end{eqnarray*}
which implies
\begin{eqnarray*}
C_N(\Omega )\leq
\frac{2}{1-\Vert P\Vert }\leq \frac{\pi }{\alpha }.
\end{eqnarray*}
\noindent{\it Remark.} For $\varepsilon>0$ we consider the Banach space
$$X_{\varepsilon}:=\{f\in C(\partial\Omega)\,;
\|f\|_{\varepsilon}<+\infty\},$$
where $\|f\|_{\varepsilon}:=\sup_{z\in\partial\Omega}\big(
(1\!+\!|z|)^\varepsilon
|f(z)|\big)$. We set
\begin{eqnarray*}
\gamma_{\varepsilon}:= \frac{1}{\pi}\,\sup_{z\in \partial\Omega }\int_{\partial \Omega_{z}}\frac{(1+|z|)^\varepsilon }{(1+|\sigma|)^\varepsilon
}\,d\arg(\sigma\!-\!z ).
\end{eqnarray*}
When $\gamma_{\varepsilon}<+\infty$, using the positivity of $P$,
it is easily seen that
$P$ acts from $X_{\varepsilon}$ into itself and that the
corresponding induced operator norm is $\gamma_{\varepsilon}$.
It can be seen that
$\limsup_{\varepsilon\to 0} \gamma_{\varepsilon}\leq1-2\alpha/\pi$.
Therefore $I\!+\!P$ is invertible in $X_{\varepsilon}$
for $\varepsilon$ small enough.
In particular, if $r$ is a rational function bounded in
$\Omega $ and satisfying $r(\infty )=0$, then the corresponding $g$ in
(\ref{cne}) belongs to some space $X_{\varepsilon}$ and thus it satisfies an
estimate of the form $|g(z)|\leq C(1\!+\! |z|)^{-\varepsilon }$.

\medskip
b) We consider now a rational function $R$ with matrix values and satisfying $\|R(z)\|\leq1$ in $\Omega$. The corresponding $G$ in the matrix-valued version of (\ref{cne}) satisfies
$\|G\|_{L^\infty(\partial\Omega)}\leq\pi/\alpha$.

We assume initially that $R(\infty )=0$. From the previous remark
$\|G(\sigma)\|\leq C(1\!+\! |\sigma|)^{-\varepsilon }$; this will
insure normal convergence in the following integrals. Then
(\ref{cne}) implies
\begin{eqnarray*}
R(z)=\frac{1}{4\pi}\int_{\partial \Omega }G(\sigma )\,\frac{\nu}{\sigma -z}\,ds,
\quad 0=\frac{1}{4\pi}\int_{\partial \Omega }G(\sigma )\,
\frac{\bar\nu }{\bar\sigma -\bar z}\,ds.
\end{eqnarray*}
That justifies to replace $z$ and $\bar z$ by $A$ and resp. $A^*$, if the
operator $A$
satisfies $\overline W(A)\subset \Omega $. We get
\begin{eqnarray*}
R(A)=\frac{1}{2}\int_{\partial \Omega }G(\sigma )\otimes\mu (\sigma
,A)\,ds.
\end{eqnarray*}
Note that this equality is still true for any constant $R$ (with the
corresponding constant $G$). Therefore it remains valid without the
restriction $R(\infty )=0$. Then the theorem follows from
\begin{eqnarray*}
\Vert R(A)\Vert \leq \tfrac{1}{2}\,\Vert G\Vert_{L^\infty(\partial  \Omega) }
\ \Vert \int_{\partial \Omega } \mu (\sigma ,A)\,ds\Vert=\frac{\pi -\alpha }{\pi }\,
 \Vert G\Vert_{L^\infty(\partial  \Omega) }\leq \frac{\pi -\alpha }{\alpha }.
\end{eqnarray*}
\end{proof}

Using Theorem \ref{unbdeux} and the general bound from
\cite{crac}, we obtain
 \begin{ccorollary}
We assume that the convex domain $\Omega \neq \C$ contains a
sector of angle $2\alpha $, with $0<\alpha \leq \pi /2$. Then
\begin{eqnarray*}
C_{cb}(\Omega)\leq\min\Big(  \frac{\pi -\alpha }{\alpha },
1+\frac{2}{\pi}\int_\alpha ^{\pi /2} \frac{\pi -x+\sin x }{\sin
x}\,dx, 57 \Big).
\end{eqnarray*}
 \end{ccorollary}

\noindent{\bf Remark 4.4.} (comparing the bounds) The first bound
is the best one if $0.346\, \pi\leq\alpha\leq 0.5\,\pi$ and the
second if  $4\,10^{-13}\pi\leq\alpha\leq0.345\, \pi$. The last has
to be used only if $0<\alpha\leq 4\,10^{-13}\pi$.
\medskip

The result \emph{a fortiori} holds if $\Omega $ is a sector of angle
$2\alpha $,
but in this case it is possible to obtain more precise estimates.

\begin{theorem}
\label{th:secto}
We assume that the convex domain is a sector $S_\alpha $ of
angle
$2\alpha $, with $0<\alpha \leq \pi /2$. Then we have the following estimates~:
\begin{eqnarray*}
C_N(S_\alpha)  \leq 2-\frac{2}{\pi}\log\tan\big(\frac{\alpha \,\pi}
{4(\pi \!-\! \alpha)}\big),\quad
\quad C_N(S_\alpha )=\tfrac{2}{\pi }\log\tfrac{1}{\alpha }+O(1)\\
\hbox{\rm and}\quad C_{cb}(S_\alpha  )\leq
\frac{\pi -\alpha }{\pi }\ C_N(S_\alpha ).\hskip 3cm
\end{eqnarray*}
\end{theorem}
\begin{proof} Without loss of generality we can assume that
$$S_\alpha  =\{z\in\C\,; z\neq 0 \mbox{ and } |\arg z|<\alpha \}.$$
The proof is based on the
following
relations
\begin{eqnarray*}
g(e^{t+i\alpha })= 2\,r(e^{t+i\alpha })-\int_{\R}r(e^{s-i\alpha
})\,q(t\!-\!
s)\,ds,\hskip 2cm\\ g(e^{t-i\alpha })= 2\,r(e^{t-i\alpha
})-\int_{\R}r(e^{s+i\alpha
})\,q(t\!-\!s)\,ds,\hskip2cm\\
\hbox{where}\quad q(t)=\frac{i}{\pi }\,\frac{d}{dt}
\log\frac{e^{2i\alpha
\nu}-e^{-t\nu}}{1+e^{-t\nu}},\qquad  \nu=\frac{\pi }{2(\pi \!-\! \alpha )}.
\end{eqnarray*}
We refer to
\cite{crdelog}
for all details of these computations.
We deduce
\begin{eqnarray*}
\Vert g\Vert _{L^\infty (S_\alpha )}\leq \big(2+\Vert q\Vert
_{L^1(\R)}\big)\,
\Vert r\Vert _{L^\infty (S_\alpha )} ;
\end{eqnarray*}
therefore
\begin{eqnarray*}
C_N(S_\alpha)  \leq 2+\Vert q\Vert _{L^1(\R)}=
2-\frac{2}{\pi}\log\tan\big(\frac{\alpha
\,\pi} {4(\pi \!-\! \alpha)}\big).
\end{eqnarray*}
For the function
$r(z)=\frac{1-z^{\pi /2\alpha} }{1+z^{\pi /2\alpha}}$ we have
$\Vert r\Vert _{L^\infty }=1$. Therefore
\begin{eqnarray*}
C_N(S_\alpha)\geq \Vert g\Vert _{L^\infty } \geq |g(e^{i\alpha })|
&=&\Big|   \int_{\R}r(e^{s-i\alpha })\,q(-\!s)\,ds\Big|+O(1)\\
&=&\Big|   \int_{\R}\frac{1\!+\! ie^{s\pi /2\alpha }}{1\!-\! ie^{s\pi
/2\alpha
}}\,q(-\!s)\,ds\Big|+O(1)\\
&=&\Big|   \int_{\R}{\rm sign}(s)\,q(-\!s)\,ds\Big|+O(1)\\
&=&\frac{2}{\pi }\,\log\frac{1}{\alpha }+O(1).
\end{eqnarray*}
\end{proof}

\noindent{\bf Remark 4.5.} We know from \cite{crde} or \cite{crac} that
$C_{cb}(S_\alpha )$ is uniformly bounded, while we have $\lim_{\alpha \to
0} C_N(S_\alpha )=+\infty $. The C. Neumann approach seems to be inappropriate
when estimating $C_{cb}(\Omega)$ for flat domains.

\subsection{Dilation theorems}\label{subs}
In this subsection we assume the invertibility of $I+P$ , but the
domain $\Omega $ may be unbounded. We obtain the following dilation
result.
\begin{theorem}\label{thm:berger} We assume that the convex domain
$\Omega $ is such that $I+P$ is an isomorphism of $C(\partial \Omega
)$
and that the operator $A\in \Lh$ satisfies
$\overline{W(A)}\subset\Omega$. Then there exists a larger Hilbert space
$K$ containing $H$, and a normal operator $N$ acting on $K$ with
spectrum $\sigma (N) \subset \partial \Omega$, such that, for all
rational functions $r$ bounded in $\Omega $,
\begin{eqnarray*}
r(A) = P_{H}g(N)|_H .
\end{eqnarray*}
Here $P_H$ is the orthogonal projection from $K$ onto $H$ and
$g=2(I\!+\! P)^{-1}r$.
\end{theorem}
\begin{proof}
It follows from (\ref{cne}) that
\begin{eqnarray*}
r(A) = \int_{\partial \Omega} g(\sigma )\, \mu(\sigma ,A)\,ds .
\end{eqnarray*}
The result follows from the
Naimark's dilation theorem \cite[page 50]{book} which shows the existence
of a spectral measure
$E$ dilating the regular positive measure $\mu (\sigma ,A)\,ds$.
\end{proof}

\noindent{\bf Remark 4.6.} If $\Omega = D$ is the unit disk and
$r(z)=z^n$, with $n\geq 1$, we have (see Remark 4.1) $g=2(I\!+\!
P)^{-1}r=2\,r$. In this case, the above theorem reduces to the
2-dilation theorem of
Berger \cite{best} which states that every $A \in \Lh$ with $W(A)\subset
\overline{D}$ satisfies
\begin{eqnarray*}
A^n=2\,P_{H}U^n|_H, \qquad \forall \,n\geq 1,
\end{eqnarray*}
for a suitable unitary operator $U$ acting on $K$.

\medskip

 \begin{ccorollary}[T.\,Kato]\label{cor:44}
Let $f$ be a rational function such
that $f(\infty
)=\infty $ and the set $\Omega := \{z
 \in \C :|f(z)| < 1\}$ is convex. Let $A\in \Lh $ be
a linear operator such that its numerical range satisfies
 $W(A)\subset \overline{\Omega}$. Then we have $W(f(A))
\subset\overline{D}$.
 \end{ccorollary}
 \begin{proof}
 We set $\rho (\bar z)=\overline{1/f(z)}$. Notice that $\rho $ is
 antiholomorphic out of $\Omega $ and $\rho (\bar\sigma )=f(\sigma )$ on
 $\partial \Omega $.
 We deduce from the Cauchy formula that, for all $k\geq 1$,
 and all
 $z\in \Omega $,
 \begin{eqnarray*}
 (f(z))^k&=&\int_{\partial \Omega }f^k(\sigma )\,\mu (\sigma
 ,z)\,ds+\frac{1}{2\pi i}\int_{\partial \Omega }f^k(\sigma
 )\frac{d\bar\sigma }{\bar\sigma -\bar z}\\
 &=&\int_{\partial \Omega }f^k(\sigma )\,\mu (\sigma
 ,z)\,ds+\frac{1}{2\pi i}\int_{\partial \Omega }\rho ^k(\bar\sigma
 )\frac{d\bar\sigma }{\bar\sigma -\bar z}\\
&=&\int_{\partial \Omega
 }f^k(\sigma )\,\mu (\sigma
 ,z)\,ds.
 \end{eqnarray*}
 Therefore $g=2\,f^k$ is the solution of the Neumann problem (\ref{cne})
 for the data $r=f^k$, i.e. $f^k=(I\!+\! P)^{-1}f^k$ (or equivalently
 $Pf^k=0$).
  From the previous theorem we deduce that
 \begin{eqnarray*}
 (f(A))^k=2\,P_H\,U^k|_H,\qquad \forall \,k\geq 1,
 \end{eqnarray*}
 where $U=f(N)$ is a normal operator with spectrum $\sigma (U)\subset
 f(\partial\Omega )\subset \partial D$. Therefore $U$ is a unitary
 operator and $U$ is a unitary 2-dilation of
 $f(A)$. Then it follows from \cite{best} that $W(f(A))
 \subset\overline{D}$.
 \end{proof}

 \noindent{\bf Remark 4.7.}  A domain $\Omega $ as in Corollary \ref{cor:44}
is called a \emph{convex lemniscate}.
The spectral mapping theorem stated in the above Corollary,
which is  due to T.\,Kato \cite{kato}, was proved here using a
different method. The key
 point in our proof is that $Pf^k=0$ for all $k\geq 1$
 (compare with \cite[Theorem 2.1]{eks}).

\smallskip

 \noindent{\bf Remark 4.8.} After a first version of this paper
has been completed,
 the recent article \cite{PuSa} was brought
to the authors' attention by John
 McCarthy and Mihai Putinar. There is some overlapping between
 Subsection \ref{subs} and some of the results of \cite{PuSa}.
\section{The similarity approach}\label{sect:5}
\subsection{The case of the disk}\label{subsect:51}

We first look to the case where $\Omega $ is the unit disk $D=\{z\in  \C\,; |z|<1\}$.
\begin{theorem}
\label{thdisk}
In the disk case we have $ C_{cb} (D)=C(D,2)=2$.
\end{theorem}
\begin{proof} Let $A$ be an operator with $\overline{W(A)}\subset D$. From the Berger thorem
(see Remark 4.6) we know that $A$ admits a 2-unitary dilation. Then a result
due to Okubo and Ando \cite{okan} states the existence of an invertible
operator $S$ such that~:
\begin{eqnarray}\label{oa}
\Vert S\Vert \,\Vert  S^{-1}\Vert
\leq 2 \mbox{ \it  and } \Vert S\,A\,S^{-1}\Vert \leq 1.
\end{eqnarray}
 Using the von Neumann inequality for the contraction $S\,A\,S^{-1}$
we deduce that
$
\|R(S\,A\,S^{-1})\|\leq 1,
$
for every matrix-valued rational function $R$ with $\|R(z)\|\leq1$ in $D$.
We deduce $C_{cb}(D)\leq2$ from
the inequality
$$\|R(A)\|\leq\|S^{-1}\|\|R(S\,A\,S^{-1})\|\|S\|.$$

 With the choice $r(z)=z$ and $A=\begin{pmatrix}  0 &2\\0 &0\end{pmatrix}$ we see that the bound 2 is attained.
\end{proof}

\subsection{The case of $2\times 2$ matrices}\label{subsect:52}

We turn now to the case of $2\times 2$ matrices. In this case
it is known that
the numerical range is an ellipse whose foci are the eigenvalues.
This ellipse is  degenerate only if the matrix is normal.
\begin{theorem}
\label{mat2}
Let $A$ be a $2\times 2$ (non normal) matrix. Then there exists a
conformal function $a$ mapping the ellipse $W(A)$ onto the unit disk
$D$ and
an invertible matrix $S$ such that
$\Vert S\Vert
\,\Vert S^{-1}\Vert \leq 2$ and $\Vert S\,a(A)\,S^{-1}\Vert \leq 1$.
\end{theorem}
\begin{proof}
If the eigenvalues of $A$ are equal, then $W(A)$ is a disk and
the result follows from (\ref{oa}). Suppose now that $A$ has
distinct eigenvalues. Since any
matrix is unitary similar to an upper triangular matrix, we can
assume that $A$ is upper triangular. Furthermore, it is clear that
if the theorem holds for a matrix $A$, then it also holds for $\lambda
A+\beta I$ for any $\lambda \neq 0$ and $\beta \in \C$. Finally, we
only have to look at matrices of the form
$A=\begin{pmatrix} 1 &\gamma \\ 0 &-1 \\ \end{pmatrix}$ and we can
furthermore assume that $\gamma >0$. Then $W(A)$ is the ellipse
of foci $1$, $-1$ and minor axis $\gamma $ (see \cite{gura} for
instance). If $\rho>1 $ is chosen such that $\gamma =\rho -1/\rho
$, then the major axis is $\rho +1/\rho$.
Also (see
\cite{he}), the function
$$
a(z)=\frac{2z}{\rho }\,\exp\big(- \sum_{n\geq 1}\, \frac{(-1)^{n+1}}{n}
\,\frac{2\,t_{2n}(z)}{1+\rho ^{4n}} \big),
$$
where $t_n$ denotes the $n^{th}$ Chebyshev polynomial, is the Riemann
conformal function mapping (the interior of) $W(A)$ onto $D$. Notice that
$a(1)=-a(-1)$; thus $a(A)=a(1)A$.

We choose now
$$
S=\begin{pmatrix} 1+a(1)^2&a(1)^2\rho -1/\rho   \\ 0
&a(1)(\rho+1/\rho )\end{pmatrix};
$$
then we have
$$
B:=Sa(A)S^{-1}=\begin{pmatrix} a(1)  &1\!-\! a(1)^2  \\ 0 &-a(1)
\end{pmatrix}.
$$
It is easy to verify that $\Vert B\Vert =1$. Some simple computations
show that
the quantity
$\Vert S\Vert\,\Vert S^{-1}\Vert $ is the largest root of the equation
\begin{eqnarray*}
\label{vi}
X^2-\frac{1+\rho ^2a(1)^2}{\rho\, a(1)}\,X+1=0,
\end{eqnarray*}
which is $X=\rho \,a(1)= 2\,
\exp\big(- \sum_{n\geq 1}\, \frac{(-1)^{n+1}}{n}
\,\frac{2}{1+\rho ^{4n}} \big).$ It is easy to verify that $1<X<2$.
\end{proof}

\begin{ccorollary}
\label{thdeux}
In the two dimensional case we have
$$ \sup_\Omega C_{cb}(\Omega, 2)=
\sup_\Omega C(\Omega, 2)=2.$$
\end{ccorollary}
\begin{proof} We know from the previous subsection that $C(D,2)= 2$, thus
it is
sufficient to prove, for any $2\times 2$ matrix $A$, that $
C_{cb}(W(A), 2)\leq
2$. Let
$R$ be a rational function bounded in $W(A)$\,; we set
$Q(z)=R(a^{-1}(z))$ and
$B=S\,a(A)\,S^{-1}$. We have
\begin{eqnarray*}
  R(A)=Q(a(A))=S^{-1}Q(B)S.
\end{eqnarray*}
Using the von Neumann inequality for the
contraction $B$, we obtain
\begin{eqnarray*}
\Vert R(A)\Vert \leq 2\,\Vert Q(B)\Vert \leq 2\,\sup_{z\in D}\Vert Q(z)
\Vert = 2\,\sup_{\zeta \in W(A) }\Vert R(\zeta  )\Vert ;
\end{eqnarray*}
thus  $C_{cb}(W(A), 2)\leq 2$.
\end{proof}

\noindent{\bf Remark 5.1.}
The equality
$C(\Omega ,2)=C_{cb}(\Omega,2)$ follows also from \cite{paul3}.
\subsection{The case of ellipses}\label{subsect:53}

\begin{theorem}\label{thm:5.4}
Let $\mathcal{E}$ be (the interior of) an ellipse in $\C$ with foci
$\mu_1$ and
$\mu_2$ and let $\gamma$ be the length of the minor axis. Let $A$
be an operator in $\Lh$. Then
$W(A) \subset \overline{\mathcal E}$ if and only if there exits an
isometry $V$ from $H$ into $H\oplus H$ such that
$$A = V^*(E\otimes I) V, \quad \mbox{ where } E = \begin{pmatrix}  \mu_1
&\gamma
     \\0 &\mu_2\end{pmatrix}.$$
\end{theorem}
\begin{proof}
The \emph{if} part is easily verified. For the \emph{only if} part we can assume without loss of generality that
$\mathcal{E}=\{x\!+\!iy\,;\cos^2\theta\,x^2+y^2<1\}$ for some real $\theta$, i.e. $\mu _1=-\mu
_2=\tan \theta$ and $\gamma =2$. Then we write
\begin{eqnarray*}
M=\tfrac{1}{2}(A+A^*),\quad N=\tfrac{1}{2i}(A-A^*);\quad{\rm thus  }\quad A=M+iN.
\end{eqnarray*}
It is easily seen that
\begin{eqnarray*}
W(A)\subset \overline{\mathcal E}\quad \Longleftrightarrow\quad
W(\cos\theta M+iN)
\subset \overline D.
\end{eqnarray*}
According to a result of Ando \cite{ando}, there exist a unitary
operator $U$ and a self-adjoint operator $B$ such that
\begin{eqnarray*}
\cos \theta  M+iN&=&2\,\sin B\ U\,\cos B.
\end{eqnarray*}
Therefore
\begin{eqnarray*} \cos \theta  M&=&\sin B\ U\,\cos
B+\cos B\ U^*\,\sin B,\\
  iN&=&\sin B\ U\,\cos B-\cos B\ U^*\,\sin B.
\end{eqnarray*}
We deduce
\begin{eqnarray*}
M+iN=\frac{1}{\cos\theta }\,W^*
\begin{pmatrix}  0 &1\!+\! \cos\theta
     \\1\!-\! \cos\theta  &0\end{pmatrix} W,\quad \hbox{with}\quad
W= \begin{pmatrix}  \sin B\\U\cos B \end{pmatrix}.
\end{eqnarray*}
It is clear that $W$ is an isometry. Now we remark that
\begin{eqnarray*}
\begin{pmatrix}  0 &1\!+\! \cos\theta
     \\1\!-\! \cos\theta  &0\end{pmatrix}=
\begin{pmatrix}  \cos\frac{\theta }{2} &-\sin\frac{\theta }{2}
     \\ \sin\frac{\theta }{2}  &\cos\frac{\theta }{2}\end{pmatrix}
\begin{pmatrix}  \sin\theta  &2\cos\theta
     \\0  &-\sin\theta \end{pmatrix}
\begin{pmatrix}  \cos\frac{\theta }{2} &\sin\frac{\theta }{2}
     \\-\sin\frac{\theta }{2}  &\cos\frac{\theta }{2}\end{pmatrix}.
\end{eqnarray*}
The theorem follows by taking
\begin{eqnarray*}
V=\begin{pmatrix}  \cos\frac{\theta }{2} &\sin\frac{\theta }{2}
     \\-\sin\frac{\theta }{2}  &\cos\frac{\theta }{2}\end{pmatrix}
   W.
\end{eqnarray*}

\end{proof}

\noindent{\bf Remark 5.2.}  Theorem \ref{thm:5.4}
is useful for estimates involving
polynomials of degree one with matrix coefficients. Noticing that
for such a polynomial $P$ we have
\begin{eqnarray*}
P(A) = V^*(P(E)\otimes I_H)V ,
\end{eqnarray*}
we deduce
\begin{eqnarray*}
\|P(A)\| \leq 2\,\sup_{z\in \mathcal{E}}\|P(z)\|,
\end{eqnarray*}
for all
polynomials $P$ of degree one with matrix coefficients.

\end{document}